\documentclass{amsart}
\usepackage[english]{babel}
\usepackage{latexsym,amsfonts,amssymb,enumerate,graphics,enumerate,amsthm}

\title[Iteration of order preserving subhomogeneous maps]
{Iteration of order preserving subhomogeneous maps on a cone}

\author[M. Akian]{Marianne Akian}
\address{INRIA, Domaine de Voluceau-Rocquencourt BP 105,
78153 Le Chesnay Cedex, France,}
\email{marianne.akian@inria.fr}
\author[S. Gaubert]{St\'ephane Gaubert}
\address{INRIA, Domaine de Voluceau-Rocquencourt BP 105,
78153 Le Chesnay Cedex, France,}
\email{stephane.gaubert@inria.fr}
\author[B. Lemmens]{Bas Lemmens}
\address{Mathematics Institute, University of Warwick, Coventry CV4 7AL, U.K.}
\thanks{The third author was supported by a TALENT-Fellowship of the
Netherlands Organization for Scientific Research (NWO)}
\email{lemmens@maths.warwick.ac.uk}
\author[R.D. Nussbaum]{Roger Nussbaum}
\address{Department of Mathematics, Hill Center, Rutgers University
New Brunswick, NJ, 08903, U.S.A.}
\thanks{The fourth author was partially supported by NSF DMS-0070829 and by
NSF INT-000152}
\email{nussbaum@math.rutgers.edu}

\date{September 2004}

\let\epsilon=\varepsilon
\newtheorem{theorem}{Theorem}[section]
\newtheorem{definition}[theorem]{Definition}
\newtheorem{lemma}[theorem]{Lemma}
\newtheorem{corollary}[theorem]{Corollary}

\begin{document}
\begin{abstract}
We investigate the iterative behaviour of continuous  order
preserving subhomogeneous maps $f\colon K\to K$, where $K$ is a
polyhedral cone in a finite dimensional vector space. We show that
each bounded orbit of $f$ converges to a periodic orbit and,
moreover, the period of each periodic point of $f$ is bounded by
\[ \beta_N = \max_{q+r+s=N}\frac{N!}{q!r!s!}=
\frac{N!}{\lfloor\frac{N}{3}\rfloor!\lfloor\frac{N+1}{3}\rfloor!
\lfloor\frac{N+2}{3}\rfloor!}\sim \frac{3^{N+1}\sqrt{3}}{2\pi N},
\]
where $N$ is the number of facets of the polyhedral cone. By
constructing examples on the standard positive cone in
$\mathbb{R}^n$, we show that the upper bound is asymptotically
sharp.

These results are an extension of work by Lemmens and
Scheutzow concerning periodic orbits in the interior of the standard 
positive cone in $\mathbb{R}^n$.
\end{abstract}

\maketitle

\section{Introduction}\label{sec:1}
Let $K$ be a polyhedral cone in a finite dimensional real vector
space $X$ and $f\colon K\to K$ be a continuous map. A 
basic problem in the theory of discrete dynamical systems is to
describe qualitatively the asymptotic behaviour of the orbits
$\{f^k(x)\colon k=0,1,2,\ldots\,\}$ for each initial point $x\in
K$, as  $k\to\infty$. In this paper we investigate this problem
for continuous maps $f\colon K \to K$ that are, in addition, order
preserving and subhomogeneous. In particular, we prove in 
Theorem \ref{thm:2.1} that each
bounded orbit of $f$ converges to a periodic orbit and that the
period of each periodic point of $f$ is bounded by
\begin{equation}
\label{eq:0} \beta_N = \max_{q+r+s=N}\frac{N!}{q!r!s!}=
\frac{N!}{\lfloor\frac{N}{3}\rfloor!\lfloor\frac{N+1}{3}\rfloor!
\lfloor\frac{N+2}{3}\rfloor!},
\end{equation}
where $N$ is the number of facets of the polyhedral cone $K$. Here
$\lfloor a\rfloor$ denotes the greatest integer not exceeding $a$.
As a second result we show in Theorem \ref{thm:2.2} that the upper bound is
asymptotically sharp in case the polyhedral cone is the standard
positive cone in $\mathbb{R}^n$ given by 
$\mathbb{R}^n_+=\{x\in\mathbb{R}^n \colon x_i\geq 0 \mbox{ for
}1\leq i\leq n\}$.

Order preserving subhomogeneous maps have been studied intensively in 
nonlinear Perron-Frobenius theory.  
They arise in various fields, such as optimal control and game theory 
\cite{AG,NS,RS}, idempotent analysis \cite{KoM,MS},
the analysis of monotone dynamical systems \cite{H,J,KN,KR,Ta1,Ta2}, and 
discrete event systems \cite{BCO,Gu1,Gu2}. 
In this list we have quoted only a few recent works and we suggest the reader 
to consult \cite{N2,N3} for further references. 
The dynamical behaviour of these maps has been investigated in
\cite{AG,GG,H,J,KN,KR,LS2,N2,N3,N1,N4,Ta1,Ta2,We}; often under the
additional assumption that $f$ leaves the interior of $K$, denoted
$\mathrm{int}(K)$, invariant. In particular, it is known that if 
$f\colon\mathrm{int}(K)\to \mathrm{int}(K)$ is an order preserving
subhomogeneous map, then $f$ is nonexpansive with respect to
Thompson's part metric (see \cite{Bu} and \cite{N2}). 
For maps that are nonexpansive with respect
to Thompson's part metric, 
Weller \cite{We} has proved that every bounded orbit in the interior 
of the polyhedral cone $K$ converges to a periodic orbit.
Moreover, for the standard positive cone, $\mathbb{R}^n_+$, it has
been shown by Martus \cite{Ma} that if $f\colon
\mathrm{int}(\mathbb{R}^n_+)\to\mathrm{int}(\mathbb{R}^n_+)$ is
nonexpansive with respect to the part metric, then the periods of
periodic points of $f$ are bounded by $n!2^n$. The upper bound of
Martus is not sharp. In fact, Nussbaum \cite[p.525]{N1} has
conjectured that $2^n$ is the optimal upper bound; but at present
this conjecture is proved only for $n\leq 3$. The case $n=3$ is
proved by Lyons and Nussbaum in \cite{LyN1}, in which also
additional evidence supporting the conjecture is given. The
current best general estimate is $\max_k 2^k {n\choose k}$ by
Lemmens and Scheutzow \cite{LS2}. Other upper bounds have been
obtained in \cite{BW,N1,Si}. For order preserving homogeneous
maps $f\colon\mathrm{int}(\mathbb{R}^n_+)\to\mathrm{int}(\mathbb{R}^n_+)$, it
was expected that stronger estimates hold for the  periods of
periodic points. Indeed, Gunawardena and Sparrow conjectured (see
\cite{Gu1}) that ${n\choose\lfloor n/2\rfloor}$ is the optimal
upper bound. A proof of this conjecture was given by Lemmens and
Scheutzow in \cite{LS2}. We shall see that the arguments in
\cite{LS2} can be refined to show that if
$f\colon\mathrm{int}(K)\to\mathrm{int}(K)$ is an order preserving
subhomogeneous map, then the periods of periodic points of $f$ do
not exceed ${N\choose \lfloor N/2\rfloor}$, where $N$ is the
number of facets of the polyhedral cone $K$. In connection with
these results it is useful to mention that each order preserving
subhomogeneous map $f\colon\mathrm{int}(K)\to K$ is continuous and
has a continuous extension $f\colon K\to K$, which is again order
preserving and subhomogeneous (see \cite[Theorem 3.10]{BNS}).

With these results in mind the following questions are natural.
Given a polyhedral cone $K$ and a continuous order preserving
subhomogeneous map $f\colon K\to K$, does every bounded orbit of $f$ 
converge to a periodic orbit? Does there exist an a priori upper
bound for the periods of periodic points in terms of the number of
facets of $K$? If so, what is the optimal upper bound? In this
paper we answer these questions.

To conclude the introduction we outline the organisation of the paper. 
In Section 2 we state the two main results: Theorems \ref{thm:2.1}
and \ref{thm:2.2}. In Section 3 we collect some preliminary
results. Subsequently, we study in Section 4 periodic points of
order preserving subhomogeneous maps on polyhedral cones, whose orbit 
is contained in a part of the cone.
Using a result of Lemmens and Scheutzow \cite{LS2} we give an
upper bound for the possible periods of these periodic points. 
This upper bound is then
used in Section 5 to show that the period of any periodic point 
does not exceed $\beta_N$, where $\beta_N$ is given in
(\ref{eq:0}). In Section 6 we prove that each bounded orbit
converges to a periodic orbit. Combining this result with the
results in Section 5 yields the first main result, Theorem
\ref{thm:2.1}. In Section 7 we prove the second main result,
Theorem \ref{thm:2.2}.

\section{Statement of the main results}\label{sec:2}
Let $X$ be a real topological vector space. A subset $K$ of $X$ is
called a \emph{cone} if it is a convex subset of $X$ such that
$\lambda K\subset K$ for all $\lambda\geq 0$ and
$K\cap(-K)=\{0\}$. A cone $K$ in $X$ is called a \emph{closed
cone} if it is a closed subset of $X$. 
If $X$ is a finite dimensional topological vector space, then it is known that 
$X$ has exactly one Hausdorff vector space topology and it coincides with the 
standard topology. The main results of this 
paper concern closed cones in finite dimensional vector spaces. 
In that case the vector space topology will always be the standard  
topology. Many preliminary results will however be stated and
proved for more general topological vector spaces.

A closed cone $K$ in a finite dimensional vector space $X$ is said
to be a \emph{polyhedral cone} if it is the intersection of
finitely many closed half spaces, i.e., there exist linear functionals 
$\varphi_1,\ldots,\varphi_m$  such that 
$K=\{x\in X\colon \varphi_i(x)\geq 0\mbox{ for }1\leq i\leq m\}$. 
A \emph{face} of a polyhedral cone $K$ is any
set of the form $F=K\cap\{x\in X\colon \varphi(x)=0\}$, where
$\varphi\colon X\to\mathbb{R}$ is a linear functional such that
$K\subset\{x\in X\colon \varphi(x)\geq 0\}$. Note that the cone
itself is a face. The \emph{dimension} of a face $F$, denoted
$\mathrm{dim}(F)$, is the dimension of its linear span. A face $F$
is called a \emph{facet} if $\mathrm{dim}(F)=\mathrm{dim}(K)-1$.
We remark that if $K$ is a polyhedral cone with $N$ facets, then
there exist $N$ linear functionals $\psi_i\colon X\to\mathbb{R}$,
where $1\leq i\leq N$, such that
\begin{equation}
\label{eq:2.0}
K=\{x\in X\colon \psi_i(x)\geq 0\mbox{ for }1\leq i\leq N\}\cap
\mathrm{span}(K)
\end{equation}
and each linear functional $\psi_i$ defines a facet of $K$ (see
\cite[Section 8.4]{Sch}). In this paper the closed cone will often
be polyhedral and we reserve the notation $\psi_i$, where $1\leq
i\leq N$, to denote the linear functionals that define its facets.
A natural example of a polyhedral cone is the \emph{standard
positive cone} in $\mathbb{R}^n$ given by
$\mathbb{R}^n_+=\{x\in\mathbb{R}^n\colon x_i\geq 0\mbox{ for
}1\leq i\leq n\}$, which has $n$ facets.

A cone $K$ in a topological vector space $X$ induces a partial ordering 
$\leq_K$ on $X$ by $x\leq_K
y$ if $y-x \in K$.  We simply write $\leq$ if $K$ is obvious from
the context. Subsets of $X$ will always inherit the partial
ordering of $X$. If $(S,\leq)$ and $(T,\leq)$ are two partial
ordered sets, then we call a map $f\colon S\to T$ \emph{order
preserving} if $f(x)\leq f(y)$ for all $x,y\in S$ with $x\leq y$.
If $X$ is a vector space with a partial ordering $\leq$ and if
$f\colon D\to X$, where $D\subset X$, has the property that
$\lambda f(x)\leq f(\lambda x)$ for every $x\in D$ and  
$0<\lambda < 1$ satisfying $\lambda x \in D$, then $f$ is said to
be \emph{subhomogeneous}. If $\lambda f(x)=f(\lambda x)$ for every
$x\in D$ and $\lambda\geq 0$ satisfying $\lambda x\in D$, then $f$
is said to be \emph{homogeneous}. 

If $S$ is a set and $f\colon S\to S$, then a point $x\in S$ is
called a \emph{periodic point} if $f^p(x)=x$ for some integer
$p\geq 1$; the minimal such $p\geq 1$ is said to be the
\emph{period} of $x$ under $f$. The \emph{orbit} of $x\in S$ under
$f$ is given by $\mathcal{O}(x;f)=\{f^k(x)\colon
k=0,1,2,\ldots\,\}$. If $x$ is a periodic point, then
$\mathcal{O}(x;f)$ is called a \emph{periodic orbit}.

Equipped with these notions we now state the main results.
\begin{theorem}
\label{thm:2.1} Let $K$ be a polyhedral cone with $N$ facets in a
finite dimensional vector space $X$. If $f\colon K\to K$ is a
continuous order preserving subhomogeneous map and the orbit of
$x\in K$ is bounded, then there exists a periodic point $\xi$ of
$f$, with period $p$, such that $\lim_{k\to\infty} f^{kp}(x)=\xi$
and $p\leq \beta_N$, where
\begin{equation}
\label{eq:2.1}
\beta_N = \max_{q+r+s=N}\frac{N!}{q!r!s!}=
\frac{N!}{\lfloor\frac{N}{3}\rfloor!\lfloor\frac{N+1}{3}\rfloor!
\lfloor\frac{N+2}{3}\rfloor!}.
\end{equation}
\end{theorem}
To show that the upper bound $\beta_N$ is asymptotically sharp we
prove in Section \ref{sec:7} the following theorem.
\begin{theorem}
\label{thm:2.2} For every $1\leq m\leq n$, $1\leq p\leq {m\choose
\lfloor m/2\rfloor}$, and $1\leq q\leq {n\choose m}$, there exists
a continuous order preserving homogeneous map $f\colon
\mathbb{R}^n_+\to\mathbb{R}^n_+$ that has a periodic point with
period equal to the least common multiple of $p$ and $q$.
\end{theorem}
From Theorem \ref{thm:2.2} it follows that
\begin{equation}
\label{eq:2.2}
\alpha_N=\max \Big{\{}\mathrm{lcm\,}(p,q)\colon
1\leq p\leq {m\choose \lfloor m/2\rfloor}, 1\leq q\leq {N\choose m},
\mbox{ and } 1\leq m\leq N\Big{\}}
\end{equation}
is a lower bound for the maximum period of periodic points of
continuous order preserving subhomogeneous map
$f\colon\mathbb{R}^N_+\to\mathbb{R}^N_+$. We show in Section
\ref{sec:7} that $\lim_{N\to\infty} \alpha_N/\beta_N =1$. This
implies that the upper bound in Theorem \ref{thm:2.1} is
asymptotically sharp in case $K$ is the standard positive cone in
$\mathbb{R}^N$. This fact is illustrated in Table \ref{tab:1}
below. Moreover, by using Stirling's formula, it can be shown
that $\beta_N$ has the following asymptotics:
\[
\beta_N\sim \frac{3^{N+1}\sqrt{3}}{2\pi N}.
\]
\begin{table}[h]
\caption{The lower and upper bound for $1\leq N\leq 15$ \label{tab:1}}
\begin{center}
\small{
\begin{tabular}{ll}\hline
$\alpha_N$ & 1, 2, 6, 12, 30, 78, 210, 540, 1660, 4180, 11480, 34510,
 90090, 251874, 756252\\
$\beta_N$ &  1, 2, 6, 12, 30, 90, 210, 560, 1680, 4200, 11550, 34650,
 90090, 252252, 756756\\ \hline
\end{tabular}}
\end{center}
\end{table}

\section{Preliminary results}\label{sec:3}
In this section we collect several preliminary results.

\subsection{Partially ordered sets}
Partially ordered sets occur frequently in this exposition and it is useful 
to recall several basic concepts concerning them.
Let $(S,\leq)$ be a partially ordered set. We say that $a$ and
$b$ are \emph{comparable} if $a\leq b$ or $b\leq a$. A subset
$\mathcal{A}$ of $S$ is called an \emph{antichain} if no two
distinct elements in $\mathcal{A}$ are comparable. A subset
$\mathcal{C}$ of $S$ is called a \emph{chain} if every two
elements in $\mathcal{C}$ are comparable, and it is said to be a
\emph{maximal chain} if there exists no chain $\mathcal{D}\subset
S$ that properly contains $\mathcal{C}$. We have the following
basic lemma.
\begin{lemma}
\label{lem:3.1} Let $(S,\leq)$ be a partially ordered set and let
$f\colon S\to S$ be an order preserving map. If $x\in S$ is a
periodic point of $f$, then $\mathcal{O}(x;f)$ is an antichain.
\end{lemma}
\begin{proof}
Let $x\in S$ be a periodic point of $f$ with period $p$. Suppose
that $y,z\in\mathcal{O}(x;f)$ and $y\leq z$. As $\mathcal{O}(x;f)$
is a periodic orbit with period $p$, there exists $0\leq k<p$ such
that $z=f^k(y)$ and hence $y\leq f^k(y)$. Since $f^k$ is order
preserving, this implies that $f^k(y)\leq f^{2k}(y)\leq \ldots\leq
f^{kp}(y)=y$ and therefore $z\leq y$. Thus $y=z$ and hence
$\mathcal{O}(x;f)$ is an antichain.
\end{proof}

\subsection{Parts and Thompson's part metric}
Let $K$ be a cone in a topological vector space $X$. For the
analysis it is convenient to define an equivalence relation $\sim$
on $K$ by $x\sim y$ if there exist constants $0<\alpha\leq \beta$
such that $\alpha x\leq y\leq \beta x$. We write $[x]$ to denote
the equivalence class of $x$. The equivalence classes in $K$ are
called \emph{parts} (or \emph{constituents}) (see \cite{BB,Tho})
and we denote the set of all parts of $K$ by $P(K)$. We say that
$x$ \emph{dominates} $y$ if there exists $\beta>0$ such that
$y\leq \beta x$. We observe that if $x\sim x'$ and $y\sim y'$,
then $x$ dominates $y$ if and only if $x'$ dominates $y'$. This
observation allows us to define a partial ordering $\preceq$ on
the set of parts, $P(K)$, in the following manner: 
$P\preceq Q$ if $x$ dominates $y$ for some $x\in Q$ and $y\in P$.

For $x,y\in K$ we define
\begin{equation}
\label{eq:3.1}
M(y/x;K)=\inf \{\beta>0\colon y\leq \beta x\}
\end{equation}
and we put $M(y/x;K)=\infty$ if the set is empty. If $K$ is
obvious from the context, we simply write $M(y/x)$. Remark that
$M(y/x)<\infty$ if and only if $x$ dominates $y$. 
Moreover, if in addition $K$ is closed, then the infimum in (\ref{eq:3.1}) 
is attained and in that case $y\leq M(y/x)x$. 
We have the following lemma.
\begin{lemma}
\label{lem:3.2} Let $K$ be a cone in a topological vector space
$X$ and let $P$ be a part of $K$. If $\mathcal{A}$ is an antichain
in the partially ordered set $(P,\leq)$ and $f\colon P\to P$ is an
order preserving subhomogeneous map, then
\begin{equation}
\label{eq:3.3} M(f(y)/f(x))\leq M(y/x) \mbox{\quad for all
}x,y\in\mathcal{A}.
\end{equation}
Moreover, if  $f(\mathcal{A})\subset\mathcal{A}$ and each
$x\in\mathcal{A}$ is a periodic point of $f$, then
\begin{equation}
\label{eq:3.3.b} M(f(y)/f(x)) = M(y/x) \mbox{\quad for all
}x,y\in\mathcal{A}.
\end{equation}
\end{lemma}
\begin{proof}
Clearly the equations (\ref{eq:3.3}) and (\ref{eq:3.3.b}) are true
if $x=y$. So, let $x,y\in\mathcal{A}$ with  $x\neq y$. As $x$ and
$y$ belong to the same part, $M(y/x)$ is finite. Consider $\lambda
>0$ such that $\lambda>M(y/x)$. Then $y\leq \lambda x$ and, since
$\mathcal{A}$ is an antichain, we have that $\lambda> 1$. Using
the fact that $f$ is order preserving and subhomogeneous, we
deduce that $\lambda^{-1}f(y)\leq f(\lambda^{-1}y)\leq f(x)$, so
that $M(f(y)/f(x))\leq \lambda$. Since this inequality holds for
all $\lambda > M(y/x)$, inequality (\ref{eq:3.3}) follows.

To prove the second assertion we assume that $f^p(x)=x$ and
$f^q(y)=y$. By applying the previous observation iteratively we
deduce for each $k\geq 1$ that $M(f^k(y)/f^k(x))\leq
M(f(y)/f(x))\leq M(y/x)$. Now by taking $k=pq$  we find that
$M(y/x)\leq M(f(y)/f(x))\leq M(y/x)$, which completes the proof.
\end{proof}

Using the function $M(y/x)$ we define a map $d_T\colon K\times K
\to [0,\infty]$ by
\begin{equation}
\label{eq:3.4}
d_T(x,y)= \log(\max\{M(y/x),M(x/y)\})
\end{equation}
for all $(x,y)\in K\times K$, with $(x,y)\neq (0,0)$, and we
put $d_T(0,0)=0$. The function $d_T$ is called
(\emph{Thompson's}) \emph{part metric} \cite{Tho}. It is
well-known that if $K$ is a closed cone, then $d_T$ is a genuine
metric on each part of the cone, but not on the whole cone.
Indeed, $d_T(x,y)$ is finite if and only if $x\sim y$. If $K$ is
not a closed cone, then in general $d_T$ is a semi-metric on
each part. Moreover, if $K$ is a closed cone in a finite
dimensional vector space $X$ and $P$ is a part of $K$, then
$(P,d_T)$ is a complete metric space and the topology coincides with the
topology induced by the standard topology on $X$. More general
results concerning the part metric can be found in
\cite{BB,N2,N3,Tho}.

To conclude this subsection we mention the relation between the
part metric and the sup-norm on $\mathbb{R}^n$ given by 
$\|z\|_\infty=\max_i|z_i|$. Consider the standard positive cone
$\mathbb{R}^n_+$ and the part corresponding to the interior of
$\mathbb{R}^n_+$. There exists an isometry from the metric space
$(\mathrm{int}(\mathbb{R}^n_+), d_T)$ onto the metric space
$(\mathbb{R}^n, \|\cdot\|_\infty)$ (cf. \cite[Proposition 1.6]{N2}). 
Indeed, one can use the map
$L\colon\mathrm{int}(\mathbb{R}^n_+)\to\mathbb{R}^n$ given by 
\begin{equation}
\label{eq:3.5} L(x)=(\log x_1,\ldots,\log x_n) \mbox{\quad for }
x=(x_1,\ldots,x_n)\in\mathrm{int}(\mathbb{R}^n_+).
\end{equation}
The inverse of $L$ is, of course, the map
$E\colon\mathbb{R}^n\to\mathrm{int}(\mathbb{R}^n_+)$  given by 
\begin{equation}
\label{eq:3.6} E(x)=(e^{x_1},\ldots,e^{x_n})\mbox{\quad for
}x=(x_1,\ldots,x_n)\in\mathbb{R}^n.
\end{equation}
To see that the map $L$ is an isometry it is convenient to first
define a map $t\colon\mathbb{R}^n\to\mathbb{R}$ by $t(x)=\max_i
x_i$ for $x\in\mathbb{R}^n$, and subsequently to remark that
\begin{equation}
\label{eq:3.7} \|x\|_\infty=\max\{t(x),t(-x)\}.
\end{equation}
Now note that if $x,y\in\mathrm{int}(\mathbb{R}^n_+)$, then
\[
M(x/y)=\inf\{\beta\geq 0\colon x\leq\beta y\}=\max_i(x_i/y_i)
\]
and $\log(\max_i(x_i/y_i))= \max_i(\log x_i - \log y_i) =
t(L(x)-L(y))$. Thus,
\begin{equation}
\label{eq:3.8} \log M(x/y)=t(L(x)-L(y))\mbox{\quad for all 
}x,y\in\mathrm{int}(\mathbb{R}^n_+), \end{equation} so that
(\ref{eq:3.4}) and (\ref{eq:3.7}) yield
\begin{equation}
\label{eq:3.7.b} d_T(x,y)=\|L(x)-L(y)\|_\infty \mbox{\quad for all
}x,y\in\mathrm{int}(\mathbb{R}^n_+).
\end{equation}
The function $t\colon\mathbb{R}^n\to\mathbb{R}$ above appears naturally in 
the study of topical functions (see Gunawardena and Keane \cite{GK}) and also
played an important role in \cite{LS2}. It has certain properties
of a norm. For instance, $t(x+y)\leq t(x)+t(y)$ for all
$x,y\in\mathbb{R}^n$; but, $t(x)\neq t(-x)$ in general.

\subsection{Nonexpansiveness and order preserving maps on a cone}
If $(C,d)$ is a metric space, then $f\colon C\to C$ is called
\emph{nonexpansive with respect to} $d$, or, simply 
\emph{$d$-nonexpansive} if
\begin{equation}
\label{eq:3.3.1} d(f(x),f(y))\leq d(x,y)\mbox{\quad for all
}x,y\in C.
\end{equation}
The map $f$ is called a \emph{$d$-isometry} if (\ref{eq:3.3.1}) is
an equality for all $x,y \in C$. Although $d_T$ is not a proper
metric on a cone $K$ (and in general only a semi-metric on each part of $K$, 
when $K$ is not closed), we say that $f\colon K\to K$ is 
\emph{$d_T$-nonexpansive} if (\ref{eq:3.3.1}) holds for $d_T$. 
Here the inequality only makes sense if the right-hand side is finite. 
In the same way we abuse terminology for the function
$t:\mathbb{R}^n\to \mathbb{R}$ given by $t(x)=\max_i x_i$. We call
a map $f\colon S\to S$, where $S\subset \mathbb{R}^n$,
\emph{$t$-nonexpansive} if $t(f(x)-f(y))\leq t(x-y)$ for all
$x,y\in S$. The map $f$ is called a \emph{$t$-isometry} if
$t(f(x)-f(y))=t(x-y)$ for all $x,y\in S$. We have the following
lemma (cf. \cite[Proposition 1.5]{N2}), which is similar to results 
in \cite{CT}.
\begin{lemma}\label{lem:3.3}
Let $K$ be a closed cone in a topological vector space $X$. If $f\colon
K\to K$ is order preserving, then $f$ is $d_T$-nonexpansive if and only if 
$f$ is subhomogeneous.
\end{lemma}
\begin{proof}
Assume first that $f$ is subhomogeneous. If $x,y\in K$ and
\begin{equation}
\label{eq:3.3.2}
 \lambda\geq\max\{M(y/x),M(x/y)\},
\end{equation}
then $y\leq \lambda x$ and $x\leq \lambda y$, so that $x\leq
\lambda y\leq \lambda^2 x$ and therefore $\lambda\geq 1$. As $f$
is order preserving and subhomogeneous, we obtain
\[
\lambda^{-1}f(y)\leq f(\lambda^{-1}y)\leq f(x) \mbox{\quad
and\quad } \lambda^{-1}f(x)\leq f(\lambda^{-1}x)\leq f(y).
\]
This implies that $\max\{M(f(y)/f(x)),M(f(x)/f(y))\}\leq \lambda$
and hence
\[
d_T(f(x),f(y))\leq \log \lambda = d_T(x,y).
\]

Now assume that $f$ is nonexpansive with respect to $d_T$ on $K$.
Let $x\in K$ and put $y=\lambda^{-1}x$, where $\lambda\geq 1$.
Clearly $d_T(x,y)=\log \lambda$ if $x\neq 0$ and $d_T(x,y)\leq \log \lambda$ 
if $x=0$. As $f$ is nonexpansive with
respect to $d_T$, we have that
\[
\log M(f(x)/f(y))\leq d_T(f(x),f(y))\leq d_T(x,y)\leq\log\lambda,
\]
so that $f(x)\leq \lambda f(y)$. This implies that
$\lambda^{-1}f(x)\leq f(y)=f(\lambda^{-1}x)$ and hence $f$ is
subhomogeneous.
\end{proof}
Maps that are nonexpansive with respect to the part metric map
parts into parts. Indeed, we have the following lemma.
\begin{lemma}\label{lem:3.4}
If $K$ is a cone in a topological vector space $X$ and $f\colon
K\to K$ is $d_T$-nonexpansive, then $f([x])\subset
[f(x)]$ for each $x\in K$.
\end{lemma}
\begin{proof}
If $y\in f([x])$, then there exists $z\in [x]$ such that $f(z)=y$.
Since $z\sim x$ we have that $d_T(x,z)$ is finite. As $f$ is
nonexpansive with respect to $d_T$ on $[x]$, we find that
$d_T(f(x),y)$ is finite and hence $y\sim f(x)$.
\end{proof}
This lemma has the following corollary.
\begin{corollary}
\label{cor:3.3.3} If $K$ is a cone in a topological vector space
$X$ and $f\colon K\to K$ is $d_T$-nonexpansive, then the map
$F\colon P(K)\to P(K)$ given by $F(P)=[f(x)]$ for $x\in P$ is well
defined. Moreover, if $f$ is order preserving and $K$ is closed, then 
$F$ preserves the ordering $\preceq$ on $P(K)$.
\end{corollary}
\begin{proof}
To see that $F$ is well defined we let $P$ be a part of the cone $K$.
For each $x,y\in P$ we have that $f(x)\sim f(y)$, by Lemma \ref{lem:3.4}, and 
hence $[f(x)]=[f(y)]$. Thus, $F$ is well defined. If $f\colon K\to K$ is an 
order preserving $d_T$-nonexpansive map and $K$ is closed, then $f$ is   
subhomogeneous by Lemma \ref{lem:3.3}. Now let $P,Q\in P(K)$ be such that
$P\preceq Q$. If $x\in Q$ and $y\in P$, then $x$ dominates $y$
and therefore there exists $\lambda\geq 1$ such that $y\leq
\lambda x$. Since $f$ is order preserving and subhomogeneous, it
follows that $\lambda^{-1}f(y)\leq f(\lambda^{-1}y)\leq f(x)$.
Thus $f(x)$ dominates $f(y)$, so that $[f(y)]\preceq [f(x)]$. From
this we conclude that $F(P)\preceq F(Q)$, which completes the
proof.
\end{proof}

\section{Periodic orbits in a part of the polyhedral cone}
\label{sec:3.b} In \cite{LS2} Lemmens and Scheutzow  proved the
following theorem.
\begin{theorem}[\cite{LS2}]
\label{thm:3.5} If $\mathcal{A}$ is a finite antichain in
$(\mathbb{R}^n,\leq)$, where the partial ordering $\leq$ is
induced by $\mathbb{R}^n_+$, and on $\mathcal{A}$ a commutative
group of $t$-isometries acts transitively, then $\mathcal{A}$ has
at most ${n\choose \lfloor n/2\rfloor}$ elements.
\end{theorem}
We use this theorem to derive the following result.
\begin{theorem}
\label{thm:3.6} Let $K$ be a polyhedral cone with nonempty
interior in a finite dimensional vector space $X$. If $K$ has $N$
facets and $f\colon\mathrm{int}(K)\to\mathrm{int}(K)$ is order
preserving and subhomogeneous, then the periods of periodic points
of $f$ do not exceed ${N\choose\lfloor N/2\rfloor}$.
\end{theorem}
\begin{proof}
Let $\xi$ be a periodic point of $f$ with period $p$ and let
$\mathcal{A}=\mathcal{O}(\xi;f)$. From Lemma \ref{lem:3.1} it
follows that $\mathcal{A}$ is an antichain in
$(\mathrm{int}(K),\leq_K)$. Furthermore Lemma \ref{lem:3.2}
implies that
\begin{equation}\label{eq:3.9}
M(f(y)/f(x);K)=M(y/x;K)\mbox{\quad for all } x,y\in\mathcal{A}.
\end{equation}
Define $\Psi\colon X\to\mathbb{R}^N$ by
$\Psi(x)=(\psi_1(x),\ldots,\psi_N(x))$ for all $x\in X$. Here
$\psi_i\colon X\to\mathbb{R}$, with $1\leq i\leq N$, are the
linear functionals that define the facets of $K$. The map $\Psi$ is 
linear and, by (\ref{eq:2.0}), $x\in K$ if and only if
$\Psi(x)\in\mathbb{R}^N_+$. Hence
$\Psi(K)=\Psi(X)\cap\mathbb{R}^N_+$ and if $\mathbb{R}^N$ is
endowed with the partial ordering induced by $\mathbb{R}^N_+$, we
get that $x\leq_K\lambda y$ is equivalent to $\Psi(x)\leq \lambda
\Psi(y)$. It follows that
\begin{equation}
\label{eq:3.10} M(y/x;K)=
M(\Psi(y)/\Psi(x);\mathbb{R}^n_+)\mbox{\quad for all } x,y\in K.
\end{equation}
Moreover, $\Psi$ is injective, because $\Psi(x)=0$ implies that
$x\in K$ and $-x\in K$, so that $x=0$.  We also have that
$\Psi(\mathrm{int}(K))\subset\mathrm{int}(\mathbb{R}^N_+)$.
Indeed, if $y\in\mathrm{int}(K)$, then for each $z\in X$ there
exists $\epsilon>0$ such that $y-\epsilon z\in K$. This implies
that $\psi_i(y)\geq\epsilon \psi_i(z)$ for all $1\leq i\leq N$.
Since $\psi_i$ is nonzero, there exists $z\in X$ such that
$\psi_i(z)>0$. Therefore $\psi_i(y)>0$ for all $1\leq i\leq N$ and
hence $\Psi(y)\in\mathrm{int}(\mathbb{R}^N_+)$.

Let $\Psi^{-1}$ be the inverse of $\Psi$ on $\Psi(X)$. Put
$\mathcal{A}'=\Psi(\mathcal{A})$ and let
$g\colon\mathcal{A}'\to\mathcal{A}'$ be given by $g=\Psi\circ
f\circ\Psi^{-1}$. By using (\ref{eq:3.9}) and (\ref{eq:3.10}) we
find that if $u,v\in\mathcal{A}'$, $u=\Psi(x)$, and $v=\Psi(y)$,
then
\begin{equation}\label{eq:3.11}
M(v/u;\mathbb{R}^n_+)=M(y/x;K)=M(f(y)/f(x);K)=M(g(u)/g(v);\mathbb{R}^n_+).
\end{equation}

Now put $\mathcal{A}''=L(\mathcal{A}')$ and define
$h\colon\mathcal{A}''\to\mathcal{A}''$ by $h=L\circ g\circ E$,
where the maps $L$ and $E$ are given in (\ref{eq:3.5}) and
(\ref{eq:3.6}), respectively. The set $\mathcal{A}''$ is well
defined, as $\mathcal{A}'\subset
\Psi(\mathrm{int}(K))\subset\mathrm{int}(\mathbb{R}^N_+)$. It
follows from (\ref{eq:3.8}) and (\ref{eq:3.11}) that
\[
t(h(r)-h(s))=t(r-s)\mbox{\quad for all } r,s\in\mathcal{A}''.
\]
One can verify that $\mathcal{A}''$ is a periodic orbit of $h$
with period $p$; in fact, $\mathcal{A}''=\mathcal{O}(L(\Psi(\xi));h)$. 
Therefore
$G=\{h^k\colon\mathcal{A}''\to\mathcal{A}''\mid 0\leq k<p\}$ is a
commutative group of $t$-isometries that acts transitively on
$\mathcal{A}''$ and hence Theorem \ref{thm:3.5} implies that
$p=|\mathcal{A}''|\leq {N\choose \lfloor N/2\rfloor}$.
\end{proof}
We shall generalize Theorem \ref{thm:3.6} to the case where $f$
maps a part of the cone into itself; but before we do this we introduce some
definitions. Let $K$ be a polyhedral cone with $N$ facets and let 
$\psi_i\colon X\to\mathbb{R}$, with $1\leq i\leq N$, be 
the linear functionals that define the facets of $K$. 
We define for each $x\in K$ a set $I_x$ by 
\begin{equation}
\label{eq:3.b.1} I_x=\{i\in\{1,\ldots,N\}\colon \psi_i(x)>0\}. 
\end{equation}
It easy to verify that $I_y\subset I_x$ if and only if $x$
dominates $y$. Therefore $I_x=I_y$ is equivalent to $x\sim y$.
This allows us to make the following definition.
\begin{definition}\label{def:3.0}
If $K$ is a polyhedral cone with $N$ facets in a finite
dimensional vector space $X$, then for each part $P\in P(K)$ we
define $I(P)=I_x$, where $x\in P$.
\end{definition}
The same observation shows that the map $I\colon P(K)\to 2^{[N]}$
given by $P\mapsto I(P)$ is injective. Here $2^{[N]}$ denotes the set of
all subsets of $\{1,\ldots,N\}$. In particular, this implies that
there are at most $2^N$ parts in $K$. Moreover, $I(P)\subset I(Q)$
if and only if $P\preceq Q$, and hence $I\colon
(P(K),\preceq)\to(2^{[N]},\subset)$ and its inverse $I^{-1}$ are both
order preserving.
\begin{corollary}
\label{cor:3.7} Let $K$ be a polyhedral cone in a finite
dimensional vector space $X$. If $P$ is a part of $K$ and $f\colon
P\to P$ is order preserving and subhomogeneous, then the periods
of periodic points of $f$ do not exceed ${m\choose \lfloor
m/2\rfloor}$, where $m=|I(P)|$.
\end{corollary}
\begin{proof}
Let $K$ be given by $\{x\in X\colon \psi_i(x)\geq 0\mbox{ for }
1\leq i\leq N\}\cap \mathrm{span}(K)$, where each $\psi_i$ is a
linear functional that defines a facet of $K$. Put $J=I(P)$ and
let $J'$ denotes its complement. Define
\[
Y=\{x\in X\colon \psi_j(x)=0 \mbox{ for all }j\in
J'\}\cap\mathrm{span}(K).
\]
Remark that $Y$ is a linear subspace of $X$. Now let $C=K\cap Y$.
We observe that $C=\{y\in Y\colon \psi_j(y)\geq 0\mbox{ for all
} j\in J\}$ and hence $C$ is a polyhedral cone with at most $|J|$
facets in the vector space $Y$. Since
\begin{eqnarray*}
P & = & \{x\in X\colon \psi_j(x)>0\mbox{ for }j\in J\mbox{ and
}\psi_j(x)=0
\mbox{ for }j\in J'\}\cap\mathrm{span}(K) \\
  & = & \{y\in Y\colon \psi_j(y)>0\mbox{ for } j\in J\}, 
\end{eqnarray*}
we have that $P$ is the interior of $C$ in $Y$. Thus we can apply
Theorem \ref{thm:3.6} to conclude that the periods of periodic
points of $f\colon P\to P$ do not exceed ${q\choose \lfloor
q/2\rfloor}$, where $q$ is the number of facets of $C$. Since
$q\leq |J|=|I(P)|=m$, we find that ${q\choose \lfloor
q/2\rfloor}\leq {m\choose \lfloor m/2\rfloor}$ and this completes
the proof.
\end{proof} Theorem \ref{thm:3.6} and Corollary \ref{cor:3.7}
generalize Theorem 5.2 in \cite{LS2} by allowing subhomogeneous
maps rather than homogeneous maps and allowing general polyhedral
cones. To conclude this section we mention one other consequence
of Theorem \ref{thm:3.6}, which refines another result in
\cite{LS2}. It concerns order preserving sup-norm nonexpansive
maps. Recall that a map $f\colon\mathbb{R}^n\to\mathbb{R}^n$ is
sup-norm nonexpansive if $\|f(x)-f(y)\|_\infty\leq \|x-y\|_\infty$
for all $x,y\in\mathbb{R}^n$.
\begin{theorem}
\label{thm:3.8} If $f\colon\mathbb{R}^n\to\mathbb{R}^n$ is a
sup-norm nonexpansive map and $f$ is order preserving with respect
to the ordering induced by $\mathbb{R}^n_+$, then the periods of
periodic points of $f$ do not exceed $n\choose \lfloor
n/2\rfloor$.
\end{theorem}
\begin{proof}
Let $f\colon\mathbb{R}^n\to\mathbb{R}^n$ be an order preserving
sup-norm nonexpansive map and suppose that $\xi\in\mathbb{R}^n$ is
a periodic point of $f$ with period $p$. Define a map $h\colon
\mathrm{int}(\mathbb{R}^n_+)\to\mathrm{int}(\mathbb{R}^n_+)$ by
$h=E\circ f\circ L$, where $L$ and $E$ are respectively given in
(\ref{eq:3.5}) and (\ref{eq:3.6}). From (\ref{eq:3.7.b}) we know
that $L$ is an isometric homeomorphism between
$(\mathrm{int}(\mathbb{R}^n_+), d_T)$ and
$(\mathbb{R}^n,\|\cdot\|_\infty)$ and the inverse isometry is the
map $E$. As $f$ is sup-norm nonexpansive this implies that $h$ is
nonexpansive with respect to $d_T$. Since $f$ is order preserving,
the map $h$ also preserves the ordering induced by
$\mathbb{R}^n_+$. Therefore it follows from Lemma \ref{lem:3.3}
that $h$ is an order preserving subhomogenous map. Clearly
$E(\xi)$ is a periodic point of $h$, with period $p$, and hence we
conclude from Theorem \ref{thm:3.6} that $p$ is at most
${n\choose\lfloor n/2\rfloor}$.
\end{proof}

\section{Periods of periodic points in a polyhedral cone}\label{sec:4}

The main goal of this section is to prove that the periods of
all periodic points of order preserving subhomogenous maps on a
polyhedral cone with $N$ facets do not exceed $\beta_N$, where
$\beta_N$ is given in (\ref{eq:2.1}). But first we show the following
theorem.
\begin{theorem}
\label{thm:4.1} Let $K$ be a polyhedral cone with $N$ facets in a
finite dimensional vector space $X$. If $f\colon K\to K$ is an
order preserving subhomogenous map and $x\in K$ is a periodic
point of $f$ with period $p$, then there exist integers $q_1$ and
$q_2$ such that $p=q_1q_2$,
\[ 1\leq q_1 \leq {N\choose\max\{m,\lfloor N/2\rfloor\}},
\mbox{\quad and \quad} 1\leq q_2 \leq {m\choose\lfloor
m/2\rfloor},
\]
where $m=\min\{|I_{f^j(x)}|\colon 0\leq j<p\}$.
\end{theorem}
\begin{proof}
It follows from Lemma \ref{lem:3.3} and Corollary \ref{cor:3.3.3}
that the map $F\colon P(K)\to P(K)$ given by $F(P)=[f(x)]$ for
$x\in P$, is well defined and order preserving. Let $I\colon
P(K)\to 2^{[N]}$ be given as in Definition \ref{def:3.0}. We denote by
$I^{-1}$ the inverse of $I$ on $I(P(K))$. Define a map $G\colon
I(P(K))\to I(P(K))$ by $G=I\circ F\circ I^{-1}$. As $F$,
$I$, and $I^{-1}$ are all order preserving, the map $G$ preserves the
partial ordering $\subset$ on $I(P(K))$.

Let $x\in K$ be a periodic point of $f$, with period $p$, and let
$m=\min\{|I_{f^j(x)}|\colon 0\leq j<p\}$. Take
$z\in\mathcal{O}(x;f)$ such that $|I_z|=m$ and put $Q=[z]$. We observe that 
$F^j(Q)=[f^j(z)]$ for all $j\geq 0$ and hence $F^p(Q)=Q$. Let $k$
be the period of $Q$ under $F$. Obviously $k$ divides $p$ and
$I(Q)$ is a periodic point of $G$ with period $k$. Let
$\mathcal{A}=\mathcal{O}(I(Q);G)$. Since
\[
G^j(I(Q))=I(F^j(Q))=I([f^j(z)])=I_{f^j(z)}\mbox{\quad for all
}j\geq 0,
\]
we have that $\mathcal{A}=\{I_{f^j(z)}\colon 0\leq j<k\}$. As $G$
is order preserving, it follows from Lemma \ref{lem:3.1} that
$\mathcal{A}$ is an antichain in $(2^{[N]}, \subset)$.

A maximal chain $\mathcal{C}$ in  $(2^{[N]},\subset)$ is a sequence of
$N+1$ subsets $A_0,A_1,\ldots,A_N$ of $\{1,\ldots,N\}$ such that
$A_0\subset A_1\subset\ldots\subset A_N$ and $|A_i|=i$ for $0\leq
i\leq N$. Hence there are exactly $N!$ maximal chains. If
$A\subset\{1,\ldots,N\}$ and $|A|=s$, then there are precisely
$s!(N-s)!$ maximal chains $\mathcal{C}$ in $(2^{[N]},\subset)$ which
contain $A$. As $\mathcal{A}$ is an antichain, each maximal chain
$\mathcal{C}$ contains at most one element of $\mathcal{A}$. Since
$m=\min\{|I_{f^j(x)}|\colon 0\leq j<p\}$, we know that $|A|\geq m$
for all $A\in\mathcal{A}$.

Now for $m\leq s\leq N$, let $\nu_s$ be the number of elements of
$\mathcal{A}$ with cardinality $s$. As each maximal chain contains
at most one element of $\mathcal{A}$ and each $A\in \mathcal{A}$
with cardinality $s$ is contained in $s!(N-s)!$ maximal chains, we
find that
\[
\sum_{s=m}^N \nu_s s!(N-s)!\leq N!\mbox{\quad so that \quad }
\sum_{s=m}^N \nu_s {N\choose s}^{-1}\leq 1.
\]
Put $M(m)=\max_{m\leq s\leq N}{N\choose s}$.
It is well-known that $M(m)={N\choose m}$ if $m\geq \lfloor N/2\rfloor$, and
$M(m)={N\choose\lfloor N/2\rfloor}$ if $0\leq m\leq \lfloor N/2\rfloor$.
From this it follows that 
\begin{equation}
\label{eq:4.0}
k=|\mathcal{A}|=\sum_{s=m}^N\nu_s\leq M(m)
={N\choose\max\{m,\lfloor N/2\rfloor\}}.
\end{equation}

For $k$, $z$, and $Q$ as above, it follows from Lemma \ref{lem:3.4} that 
$f^k(Q)\subset [f^k(z)]=Q$. 
As $z$ is a periodic point of $f^k$, we can use Corollary \ref{cor:3.7} to see
that the period of $z$ under $f^k$ is less than or equal to
${m\choose \lfloor m/2\rfloor}$. Put $q_1=k$ and let $q_2$ be the
period of $z$ under $f^k$. Since $k$ divides the period $p$ of $z$
under $f$, we get that $p=kq_2=q_1q_2$, which completes the proof.
\end{proof}
We would like to remark that the arguments to derive inequality
(\ref{eq:4.0}) in the proof of Theorem \ref{thm:4.1} appear in the
study of Sperner systems and are known in combinatorics as the LYM
technique; see \cite[p. 10-11]{Bo}.

As a consequence of Theorem \ref{thm:4.1} we find that if $K$ is
a polyhedral cone with $N$ facets, then the periods of periodic
points of order preserving subhomogeneous maps $f\colon K\to
K$ are bounded by
\[
\max_{1\leq m\leq N}
{N\choose \max\{m,\lfloor N/2\rfloor\}}{m\choose \lfloor m/2\rfloor}.
\]
To see that this upper bound coincides with $\beta_N$, where
$\beta_N$ is given in (\ref{eq:2.1}), we prove the following equalities.
\begin{lemma}
\label{lem:4.2}
For each $n\geq 1$ we have that
\[
\max_{1\leq m\leq n} {n\choose \max\{m,\lfloor
n/2\rfloor\}}{m\choose \lfloor m/2\rfloor} =
\max_{q+r+s=n}\frac{n!}{q!r!s! } =
\frac{n!}{\lfloor\frac{n}{3}\rfloor!\lfloor\frac{n+1}{3}\rfloor!
\lfloor\frac{n+2}{3}\rfloor!}.
\]
\end{lemma}
\begin{proof}
We first remark that for $1\leq m\leq\lfloor n/2\rfloor$ we have that
\[
{n\choose \max\{m,\lfloor n/2\rfloor\}}{m\choose \lfloor
m/2\rfloor} \leq {n\choose \lfloor n/2\rfloor}{\lfloor
n/2\rfloor\choose \lfloor\lfloor n/2\rfloor/2\rfloor},
\]
so that
\begin{equation}
\label{eq:4.2} \max_{1\leq m\leq n} {n\choose \max\{m,\lfloor
n/2\rfloor\}}{m\choose \lfloor m/2\rfloor} = \max_{\lfloor n/2\rfloor \leq
m\leq n}{n\choose m}{m\choose \lfloor m/2\rfloor}.
\end{equation}
Further we have that
\begin{equation}
\label{eq:4.3}
{n\choose m}{m\choose \lfloor m/2\rfloor} = \frac{n!}{q!r!s!},
\end{equation}
where $q=n-m$, $r=\lfloor m/2\rfloor$, and $s=m-\lfloor m/2\rfloor$.
This implies that
\begin{equation}
\label{eq:4.4} \max_{\lfloor n/2\rfloor \leq m\leq n}{n\choose
m}{m\choose \lfloor m/2\rfloor} \leq
\max_{q+r+s=n}\frac{n!}{q!r!s!}.
\end{equation}

Let us now consider the right-hand side of (\ref{eq:4.4}).
Assume that the maximum is attained for $0\leq q^*\leq r^*\leq s^*$.
We claim that $s^*\leq q^*+1$.
Indeed, suppose by way of contradiction that  $s^* > q^*+1$.
Then $q^*!s^*!=q^*!(s^*-1)!s^*>(q^*+1)!(s^*-1)!$, so that
\[
\frac{n!}{q^*!r^*!s^*!} < \frac{n!}{(q^*+1)!r^*!(s^*-1)!},
\]
which contradicts the maximality assumption.

Since $n=q^*+r^*+s^*$ and $q^*\leq r^*\leq s^*\leq q^*+1$ we have
that $3q^*\leq n\leq 3q^*+2$ and hence $q^*=\lfloor
\frac{n}{3}\rfloor$. Furthermore, $n+1=r^*+s^*+q^*+1$ and $r^*\leq
s^*\leq q^*+1\leq r^*+1$, as $q^*\leq r^*$. This implies that
$3r^*\leq n+1\leq 3r^*+2$ and hence $r^*=\lfloor
\frac{n+1}{3}\rfloor$. Similarly, $n+2= s^*+q^*+1+r^*+1$ and
$s^*\leq q^*+1\leq r^*+1\leq s^*+1$ imply $3s^*\leq n+2\leq
3s^*+2$, so that $s^*=\lfloor\frac{n+2}{3}\rfloor$. Thus, we find
that
\begin{equation}
\label{eq:4.5}
\max_{q+r+s=n}\frac{n!}{q!r!s! } =
\frac{n!}{\lfloor\frac{n}{3}\rfloor!\lfloor\frac{n+1}{3}\rfloor!
\lfloor\frac{n+2}{3}\rfloor!}.
\end{equation}

Now put
$m=\lfloor\frac{n+1}{3}\rfloor+\lfloor\frac{n+2}{3}\rfloor$ and
compute $q$, $r$, and $s$ in the right-hand side of
(\ref{eq:4.3}). As $2\lfloor\frac{n+1}{3}\rfloor\leq m\leq
2\lfloor\frac{n+1}{3}\rfloor+1$, we find that $r=\lfloor
m/2\rfloor =\lfloor \frac{n+1}{3}\rfloor=r^*$. Moreover,
$s=m-\lfloor m/2\rfloor=\lfloor\frac{n+2}{3}\rfloor=s^*$. Since
$n=q+r+s$ we also have that $q=\lfloor \frac{n}{3}\rfloor=q^*$.
Further we remark that $m=n-q=n-\lfloor n/3\rfloor\geq 2n/3\geq
n/2\geq \lfloor n/2 \rfloor$ so that equation (\ref{eq:4.3})
implies
\[
\max_{\lfloor n/2\rfloor \leq m\leq n}{n\choose m}{m\choose\lfloor m/2\rfloor}
\geq \frac{n!}{\lfloor\frac{n}{3}\rfloor!\lfloor\frac{n+1}{3}\rfloor!
\lfloor\frac{n+2}{3}\rfloor!}.
\]
Finally we combine this inequality with (\ref{eq:4.2}), (\ref{eq:4.4}), and
(\ref{eq:4.5}) to obtain the desired result.
\end{proof}

A combination of Theorem \ref{thm:4.1} and Lemma \ref{lem:4.2}
immediately gives the following corollary.
\begin{corollary}
\label{cor:4.3}
If $K$ is a polyhedral cone with $N$ facets in a finite dimensional vector
space $X$, then the periods of periodic points of order preserving
subhomogenous maps $f\colon K\to K$ do not exceed $\beta_N$, where $\beta_N$
is given in (\ref{eq:2.1}).
\end{corollary}

\section{Asymptotic behaviour of bounded orbits}\label{sec:5}
In this section we prove Theorem \ref{thm:2.1}.
To establish this result we need to understand the asymptotic
behaviour of bounded orbits.
It is therefore natural to study the structure of the $\omega$-limit sets.
If $D$ is a metrizable topological space and $f\colon D\to D$ is a  
continuous map, then for each $x\in D$ the \emph{$\omega$-limit set} of $x$ 
under $f$ is given by
\[
\omega(x;f)=\{y\in D\colon f^{k_i}(x)\to y
\mbox{ for some sequence } (k_i)\mbox{ with }k_i\to \infty\}.
\]
It is easy to verify that each $\omega$-limit set is a (possibly empty)
closed subset of $D$ and that $f(\omega(x;f))\subset\omega(x;f)$.
Furthermore, if $\mathcal{O}(x;f)$ has a compact closure, then $\omega(x;f)$
is a nonempty compact subset of $D$ and $f(\omega(x;f))=\omega(x;f)$.
The $\omega$-limit sets also enjoy the following elementary property.
\begin{lemma}
\label{lem:5.0} Let $D$ be a metrizable topological space. If
$f\colon D\to D$ is a continuous map and $x\in D$ is such that
$\mathcal{O}(x;f)$ has a compact closure and $\omega(x;f)$ is
finite, then there exists a periodic point $\xi$ of $f$, with
period $p$, such that $\lim_{k\to\infty} f^{kp}(x)=\xi$ and
$\omega(x;f)=\mathcal{O}(\xi;f)$.
\end{lemma}
\begin{proof}
Since $f$ is continuous and $\mathcal{O}(x;f)$ has a compact closure, 
$f(\omega(x;f))=\omega(x;f)$. 
As $\omega(x;f)$ is a finite set, this implies that each $y\in\omega(x;f)$ is 
a periodic point of $f$.  
Moreover, as $\omega(x;f)$ is finite and $D$ is a metrizable topological 
space, there exist pairwise disjoint neighbourhoods $U_y$ for each 
$y\in\omega(x;f)$.
Every $y\in\omega(x;f)$ also has a neighbourhood $V_y\subset U_y$ such 
that for each $u\in V_y$ we have that $f^q(u)\in U_y$, where $q$ is the 
period of $y$ under $f$, because $f$ is continuous. 

Let $\mathrm{cl}(\mathcal{O}(x;f))$ denote the closure of 
$\mathcal{O}(x;f)$ in $D$.
Then there exists $m\geq 1$ such that for all $k\geq m$ we have that
$f^k(x)\in V_y$ for some $y\in\omega(x;f)$.
Indeed, if such an integer $m$ does not
exists, then there exists a sequence $(k_i)_i$ such that
$k_i\to\infty$ and $f^{k_i}(x)\not\in V_y$ for all
$y\in\omega(x;f)$. But $\mathrm{cl}(\mathcal{O}(x;f))$ is compact,
so that $(f^{k_i}(x))_i$ has a convergent subsequence, which has
its limit outside $\omega(x;f)$. This is obviously a
contradiction.

Now let $m \geq 1$ be such an integer. Suppose that
$f^{m}(x)\in V_z$ and let $p$ be the period of $z$. Then
$f^{m+p}(x)\in U_z$; but, as the neighbourhoods $U_y$ are
pairwise disjoint and $f^{m+p}(x)\in V_y$ for some $y\in
\omega(x;f)$, we find that $f^{m+p}(x)\in V_z$. By iterating
the argument we deduce that $f^{m+kp}(x)\in V_z$ for all
$k\geq 1$. As $z$ is the only limit point of $(f^k(x))_k$ in
$V_z$, we conclude that $(f^{m+kp}(x))_k$ converges to $z$.
This implies that $(f^{kp}(x))_k$  converges to $f^r(z)$, where
$r\equiv -m \bmod p$, because $f$ is continuous. Thus, if we take
$\xi=f^r(z)$, then $\omega(x;f)=\mathcal{O}(\xi;f)$ and this
completes the proof.
\end{proof}

To prove Theorem \ref{thm:2.1} we first show that if $f\colon K\to K$ is a 
continuous order preserving subhomogenous map on a polyhedral cone, then the 
$\omega$-limit sets of points, with a bounded orbit, are finite. A combination 
of this result with Lemma \ref{lem:5.0}, and Corollary \ref{cor:4.3} will 
yield Theorem \ref{thm:2.1}.

We shall use the following result of
Nussbaum \cite[Corollary 2]{N1}.
\begin{theorem}[\cite{N1}]
\label{thm:5.1}
Let $P$ be a part of a polyhedral cone $K$ in a finite dimensional vector
space $X$ and let $m=|I(P)|$.
If $C$ is a compact subset of $P$ and $f\colon C\to C$ is nonexpansive with
respect to $d_T$, then there exists an integer $\tau_m$, which only depends
on $m$, such that $|\omega(x;f)|\leq \tau_m$ for every $x\in C$.
\end{theorem}
The first ideas for this theorem go back to Weller \cite[Corollary
4.10]{We}, who proved a similar assertion, only without the upper
bound.

In case $K$ is the standard positive cone $\mathbb{R}^n_+$ and $P$
is the part corresponding to $\mathrm{int}(\mathbb{R}^n_+)$, we
know that the map $f\colon C\to C$, with $C\subset P$, is nonexpansive
with respect to $d_T$ if and only if the map $g\colon C'\to C'$
given by $g=L\circ f\circ E$, where $L$ and $E$ are given in
(\ref{eq:3.5}) and (\ref{eq:3.6}), is nonexpansive with respect to
the sup-norm. Using this observation it is not hard to show that
Theorem \ref{thm:5.1} is equivalent to the following assertion: if
$C$ is a compact set of $\mathbb{R}^n$ and $g\colon C\to C$ is
sup-norm nonexpansive, then there exists an integer $\tau_n$,
which only depends on $n$, such that $|\omega(x;g)|\leq \tau_n$
for all $x\in C$. It has been conjectured by Nussbaum \cite[p.
525]{N1} that the optimal choice for $\tau_n$ is $2^n$; but at
present the conjecture is proved only for $n\leq 3$ (see
\cite{LyN1}). The current best general estimate for $\tau_n$ is
$\max_k 2^k{n\choose k}$ (see \cite{LS2}).

We know by Lemma \ref{lem:3.3} that every order preserving
subhomogeneous map is nonexpansive with respect to the part
metric. Therefore Theorem \ref{thm:5.1} implies that if $f\colon
K\to K$ is an order preserving subhomogenous map and
$\mathcal{O}(x;f)\subset P$ has a compact closure in $P$, then
$\omega(x;f)$ is finite. A difficulty arises  when
$\mathcal{O}(x;f)$ is an orbit in a part $P$, but its closure is
not contained in $P$. To overcome this and other difficulties we shall use 
several technical lemmas.

Let us first recall the following old result of Freudenthal and Hurewicz 
\cite{FH}.
\begin{lemma}[\cite{FH}]
\label{lem:5.4}
Let $(C,d)$ be a compact metric space and let $D$ be a nonempty subset of
$C$.
If $f\colon D\to D$ is nonexpansive and $f$ maps $D$ onto
itself, then $f$ has a unique continuous extension
$F\colon\mathrm{cl}(D)\to\mathrm{cl}(D)$, where
$\mathrm{cl}(D)$ denotes the closure of $D$, and $F$ is an isometry
of $\mathrm{cl}(D)$ onto itself.
\end{lemma}

A combination of this lemma with Lemma \ref{lem:5.0} and Theorem \ref{thm:5.1} 
yields the following corollary. 
\begin{corollary}\label{cor:5.5.1}
Let $P$ be a part of a polyhedral cone $K$ in a finite dimensional vector 
space $X$. 
If $C$ is a compact subset of $P$ and $f\colon C\to C$ is $d_T$-nonexpansive 
map that maps $C$ onto itself, then every point $x\in C$ is a periodic point 
of $f$.
\end{corollary}
\begin{proof}
We first remark that as $C$ is a compact subset of a part of the cone, 
$(C,d_T)$ is a compact metric space. 
Since $f$ is $d_T$-nonexpansive and maps $C$ onto itself, it follows from 
Lemma \ref{lem:5.4} that $f$ is an isometry with respect to $d_T$. 

Now let $x\in C$ and remark that $\omega(x;f)$ is finite by Theorem 
\ref{thm:5.1}. As $\mathcal{O}(x;f)\subset C$, it has a compact closure. 
Therefore Lemma \ref{lem:5.0} implies that there exists a periodic point 
$\xi\in C$ of $f$, with period $p$, such that $f^{kp}(x)$ converges to 
$\xi$, as $k$ goes to infinity. 
Since $f$ is an isometry with respect to $d_T$, we find that 
\[
d_T(f^p(x),x)=d_T(f^{(k+1)p}(x),f^{kp}(x))\mbox{\quad for all } k\geq 0.
\] 
We now observe that the right-hand side of this equality 
converges to 0, as $k$ goes to infinity, and hence $d_T(f^p(x),x)=0$. Thus 
$f^p(x)=x$ and this completes the proof.
\end{proof}

The following technical lemma is stated in
considerably greater generality than is actually needed here.
Recall (see \cite[p. 41]{BNS}) that if $K$ is a closed cone in a
topological vector space $X$ and $x\in K$, we say that $K$
\emph{satisfies condition G at $x$} if for every $0<\lambda < 1$ and
every sequence $(x_k)_k$ in $K$ such that $\lim_{k\to\infty}x_k=x$, there 
exists $k^*\geq 1$ such that $\lambda x\leq x_k$ 
for all $k\geq k^*$.
We say that $K$ \emph{satisfies condition G} if it satisfies condition G at
every $x\in K$. 
If $K$ has a nonempty interior, then $K$ satisfies condition G at every 
point in its interior. 
If $K$ is a closed cone in a Hausdorff topological vector space $X$, it is
proved in \cite[Lemma 3.3]{BNS} that $K$ is a polyhedral cone in $X$
if and only if $X$ is finite dimensional and $K$ satisfies condition $G$.
\begin{lemma}
\label{lem:5.2}
Let $K$ be a closed cone in a metrizable topological vector space $X$ and let 
$D\subset K$ be such that $\lambda D\subset D$ for all $0<\lambda < 1$.
Suppose that $f\colon D\to D$ is order preserving and let $x\in D$.
If at every periodic point $\eta\in\omega(x;f)$ of $f$ condition G is 
satisfied and there exists
$\delta=\delta(\eta)>0$ such that $\lambda f^m(\eta)\leq f^m(\lambda \eta)$
for all $m\geq 1$ and $1-\delta\leq \lambda < 1$, then for every
$y\in\omega(x;f)$ and for every periodic point $\xi\in\omega(x;f)$ of $f$ 
there exists $j\geq 0$ such that $f^j(\xi)\leq y$.
Moreover, if there exists a periodic point $\xi$ in $\omega(x;f)$, then 
$\mathcal{O}(\xi;f)$ is the only periodic orbit of $f$ in $\omega(x;f)$.
\end{lemma}
Before proving this lemma we remark that if $K$ is a
polyhedral cone in a finite dimensional vector space, then $K$ satisfies
condition G at every point in $D$.
Furthermore, the condition concerning the existence of $\delta$
in Lemma \ref{lem:5.2} holds for every $y\in D$ if $f$ is subhomogeneous.
\begin{proof}[Proof of Lemma \ref{lem:5.2}]
Assume that $\xi\in\omega(x;f)$ is a periodic point of $f$ with period $p$.
By definition, there exists a sequence $(k_i)_i$
such that $f^{k_i}(x)\to\xi$, as $i\to\infty$.
By taking a subsequence we may assume that there exists  $0\leq \sigma<p$
such that $k_i\equiv \sigma \bmod p$ for all $i\geq 1$.
Take $\lambda$ with $1-\delta(\xi)\leq \lambda <1$.
As $K$ satisfies condition G at $\xi$, we have that
$\lambda \xi\leq f^{k_i}(x)$ for all sufficiently large $i$.
Suppose that $y\in\omega(x;f)$ and let $(m_i)_i$ be such that $
f^{m_i}(x)\to y$  as $i\to \infty$.
By taking a subsequence we may assume that $m_i > k_i$ for all $i\geq 1$ and
that there exists an integer $0\leq \tau<p$ such that
$m_i-k_i\equiv \tau\bmod p$ for all $i\geq 1$.
For sufficiently large $i$ we now find that 
\[
f^{m_i}(x)=f^{m_i-k_i}(f^{k_i}(x))\geq f^{m_i-k_i}(\lambda \xi) \geq
\lambda f^{m_i-k_i}(\xi)=\lambda f^\tau(\xi).
\]
Letting $i$ go infinity on the left-hand side we find that
$\lambda f^\tau(\xi)\leq y$.
Subsequently by letting $\lambda$ approach $1$ we deduce that
$f^\tau(\xi)\leq y$, which proves the first assertion.

To show the second assertion we suppose that $\xi$ and $\eta$ in 
$\omega(x;f)$ are periodic points of $f$ with period $p$ and $q$, 
respectively. 
We need to show that $\mathcal{O}(\eta;f)=\mathcal{O}(\xi;f)$.
It follows from the first assertion that there exist $0\leq \mu <p$ and 
$0\leq \nu<q$ such that $f^\mu(\xi)\leq \eta$ and $f^\nu(\eta)\leq \xi$.
Since $f$ is order preserving, it follows that $f^{\mu+k}(\xi)\leq f^k(\eta)$
and  $f^{\nu+k}(\eta)\leq f^k(\xi)$ for all $k\geq 0$.
This implies that
\[
f^{\mu+\nu}(\xi)\leq f^{\nu}(\eta)\leq \xi.
\]
By Lemma \ref{lem:3.1}, we know that $\mathcal{O}(\xi;f)$ is an antichain, so
that $\xi =f^{\mu+\nu}(\xi)$ and hence $\xi = f^{\nu}(\eta)$.
As $\eta$ and $\xi$ are both periodic points of $f$, it follows that
$\mathcal{O}(\eta;f)=\mathcal{O}(\xi;f)$ and this completes the proof.
\end{proof}

The following two lemmas tell us that we can reduce the problem to the case 
where the $\omega$-limit set is contained in a part of the cone. 
\begin{lemma}
\label{lem:5.5.2}
Let $K$ be a polyhedral cone in a finite dimensional vector space $X$. 
If $f\colon K\to K$ is $d_T$-nonexpansive, then there exists $m\geq 1$ such 
that $f^{2m}(x)\sim f^m(x)$ for all $x\in K$. 
\end{lemma}
\begin{proof}
Let $F\colon P(K)\to P(K)$ be the map in Corollary \ref{cor:3.3.3}. 
Since $P(K)$ is a finite set, we know for each $P\in P(K)$ that the sequence 
$(F^k(P))_k$ is eventually periodic, i.e., there exist $r\geq 0$ and $p\geq 1$ 
such that $F^r(P)=F^{r+kp}(P)$ for all $k\geq 0$. By the pigeonhole 
principle, we can take $r+p\leq 2^N$, where $N$ is the number of facets 
of $K$, because $|P(K)|\leq 2^N$.  Now put $m=\mathrm{lcm}\,(1,\ldots,2^N)$. 
Clearly, $r\leq m$ and $p$ divides $m$. Therefore $F^m(P)=F^{2m}(P)$ for 
each $P\in P(K)$. By taking $P=[x]$, we find that $[f^m(x)]=F^m(P)=
F^{2m}(P)=[f^{2m}(x)]$ and from this we conclude that $f^m(x)\sim f^{2m}(x)$ 
for all $x\in K$.    
\end{proof}
\begin{lemma}
\label{lem:5.5.3}
Let $K$ be a polyhedral cone in a finite dimensional vector space $X$ and let 
$g\colon K\to K$ be a continuous order preserving subhomogeneous map. 
If $x\in K$ is such that $\mathcal{O}(x;g)$ is bounded and $\mathcal{O}(x;g)$ 
is contained in a part of $K$, then $\omega(x;g)$ is contained in a part of 
$K$.  
\end{lemma}
\begin{proof}
Assume that $\mathcal{O}(x;g)$ is contained in a part $P$ of $K$.
If $P=\{0\}$, then $\mathcal{O}(x;g)=\omega(x;g)=\{0\}$ and hence the result 
is trivial in that case. Now assume that $P\neq\{0\}$, so that $I(P)$ is 
nonempty.  
We first show that there exists $c\geq 1$ such that $y\leq cx$ for all 
$y\in\mathcal{O}(x;g)$. 
As $\mathcal{O}(x;g)\subset P$, we get that $I_y=I_x=I(P)$ for all 
$y\in\mathcal{O}(x;g)$. 
This implies that $\psi_i(y)>0$ if and only if $i\in I(P)$.  
Define a number $c$ by 
\[
c=\sup\{\psi_i(y)/\psi_i(x)\colon y\in\mathcal{O}(x;g)\mbox{ and }i\in I(P)\}. 
\]
The number $c$ is finite, because $\psi_i(x)> 0$ for all $i\in I(P)$ and 
$\mathcal{O}(x;g)$ is a bounded subset of $X$. Moreover, $c\geq 1$, as 
$x\in\mathcal{O}(x;g)$ and $I(P)$ is nonempty. 
Now let $y\in\mathcal{O}(x;g)$. By definition of $c$ we have that 
$\psi_i(y-cx)\leq 0$ for all $i\in I(P)$. Since $\psi_i(y)=\psi_i(x)=0$ for 
all $i\not\in I(P)$, we deduce that $y\leq cx$. 

We remark that $\{y\in K\colon y\leq cx\}$ is a closed set that contains 
$\mathcal{O}(x;g)$ and hence it also contains $\omega(x;g)$. 
As $g$ is an order preserving subhomogeneous map and $c\geq 1$, we find that 
$g^k(y)\leq g^k(cx)\leq cg^k(x)$ for all $y\in\omega(x;g)$ and $k\geq 0$. 
The map $g$ maps $\omega(x;g)$ onto itself, because $g$ is continuous and 
$\omega(x;g)$ is bounded. Therefore $y\leq cg^k(x)$ for all 
$y\in\mathcal{O}(x;g)$  and $k\geq 0$. As the set 
$\{z\in K\colon c^{-1}y\leq z\}$ is closed, this implies that $y\leq cz$ for 
all $y,z\in\omega(x;g)$. 
Therefore $y\sim z$ for all $y,z\in\omega(x;g)$, which completes the 
proof. 
\end{proof}
Equipped with these lemmas we can now prove the following theorem.
\begin{theorem}
\label{thm:5.5}
Let $K$ be a polyhedral cone in a finite dimensional vector space $X$. If  
$f\colon K\to K$ is a continuous order preserving subhomogeneous map and 
$x\in K$ has a bounded orbit under $f$, then $\omega(x;f)$ is finite.
\end{theorem}
\begin{proof}
It follows from Lemma \ref{lem:3.4} that $f$ is $d_T$-nonexpansive. 
Let $m$ be as in Lemma \ref{lem:5.5.2} and put $g=f^m$ and $P=[g(x)]$. 
Clearly, $g(P)\subset [f^{2m}(x)]=[f^m(x)]=P$ and hence 
$\mathcal{O}(g(x);g)\subset P$. As 
$\mathcal{O}(g(x);g)\subset\mathcal{O}(x;f)$, 
the orbit $\mathcal{O}(g(x);g)$ is bounded.  
It is easy to verify that 
\[
\omega(x;f)=\omega(g(x);f)=\bigcup_{j=0}^{m-1}f^j(\omega(g(x);g)).
\]
Therefore it suffices to show that $\omega(x';g)$ is finite, whenever 
$\mathcal{O}(x';g)$ is bounded and contained in a part $P$, such that 
$g(P)\subset P$. 

So, suppose that $\mathcal{O}(x';g)$ is a bounded orbit that is contained in 
a part $P$ of $K$ and $g(P)\subset P$. It follows from Lemma \ref{lem:5.5.3} 
that $\omega(x';g)$ is included in a part of $K$, say $Q$.  
Since $\omega(x';g)$ is a bounded closed set in $Q$, we have that 
$(\omega(x';g),d_T)$ is a compact metric space. 
The map $g$ is $d_T$-nonexpansive on $Q$ and $g$ maps $\omega(x';g)$ onto 
itself. 
Therefore we can apply Corollary \ref{cor:5.5.1} and conclude that each point 
in $\omega(x';g)$ is a periodic point of $g$. As $g$ is an order preserving 
subhomogeneous maps and $K$ is a polyhedral cone, it follows from 
Lemma \ref{lem:5.2} that there is at most one periodic orbit in 
$\omega(x';g)$.  
This implies that $\omega(x';g)$ is finite and hence the proof is complete. 
\end{proof}
Knowing Theorem \ref{thm:5.5} it is now straightforward to prove 
Theorem \ref{thm:2.1}.
\begin{proof}[Proof of Theorem \ref{thm:2.1}]
Let $f\colon K\to K$ be a continuous order preserving subhomogeneous map,
where $K$ is a polyhedral cone with $N$ facets in a finite dimensional vector
space $X$.
Suppose that the orbit of $x\in K$ is bounded.
Then it follows from Theorem \ref{thm:5.5} that $\omega(x;f)$ is finite.
Therefore Lemma \ref{lem:5.0} implies that there exists a periodic point
$\xi\in K$ of $f$, with period $p$, such that $(f^{kp}(x))_k$
converges to $\xi$.
To finish the proof we remark that it follows from Corollary \ref{cor:4.3}
that $p$ is bounded by $\beta_N$, where $\beta_N$ is given in (\ref{eq:2.1}).
\end{proof}

\section{A lower bound for the maximal period}\label{sec:7}
In this section a proof of Theorem \ref{thm:2.2} is presented.
Indeed, given $1\leq m\leq n$ we construct for every $1\leq p\leq
{m\choose\lfloor m/2\rfloor}$ and $1\leq q\leq {n\choose m}$ a
continuous order preserving homogeneous map
$f\colon\mathbb{R}^n_+\to\mathbb{R}^n_+$ that has a periodic point
with period $\mathrm{lcm}\,(p,q)$. 
In the proof of Theorem \ref{thm:2.2} we use the following consequence of
an observation of Gunawardena and Sparrow (see \cite[p. 152]{Gu2}).
\begin{lemma}[\cite{Gu2}]
\label{lem:6.1} For each $1\leq p\leq {n\choose \lfloor
n/2\rfloor}$ there exists a continuous order preserving
homogeneous map $h\colon \mathbb{R}^n_+\to\mathbb{R}^n_+$ that has
a periodic point $x$, with period $p$, such that
$\mathcal{O}(x;h)\subset\mathrm{int}(\mathbb{R}^n_+)$.
\end{lemma}
Indeed, Gunawardena and Sparrow \cite{Gu2} constructed for every
$1\leq p\leq {n\choose \lfloor n/2\rfloor}$ a so called topical
map $f:\mathbb{R}^n\to\mathbb{R}^n$ that has a periodic point
$u$ with period $p$. By defining $h=\mathrm{E}\circ
f\circ\mathrm{L}$, where
$\mathrm{L}\colon\mathrm{int}(\mathbb{R}^n_+)\to\mathbb{R}^n$ and
$\mathrm{E}\colon\mathbb{R}^n\to\mathrm{int}(\mathbb{R}^n_+)$ are
respectively given in (\ref{eq:3.5}) and (\ref{eq:3.6}), we obtain
an order preserving homogeneous map
$h\colon\mathrm{int}(\mathbb{R}^n_+)\to\mathrm{int}(\mathbb{R}^n_+)$,
which has $\mathrm{E}(u)$ as a periodic point with period $p$. To
derive the conclusion of Lemma \ref{lem:6.1} we take a continuous
extension of $h$ to $\mathbb{R}^n_+$ that is order preserving
and homogeneous. Such extensions always exist (see \cite{BNS}). 
Indeed, in our case it is straightforward to find one.

The proof of Theorem \ref{thm:2.2} is quite technical.
For the reader's convenience we have therefore worked out an illustrative
example in the paragraph directly following the proof.
It may be helpful to read the two in parallel. 
Before we start the proof it useful to introduce the following notation: 
for $a,b\in\mathbb{R}$ we write $a\wedge b$ to denote $\min\{a,b\}$ and 
$a\vee b$ to denote $\max\{a,b\}$.  
\begin{proof}[Proof of Theorem \ref{thm:2.2}]
Consider a collection of $q$ distinct vectors
$\{v^1,\ldots, v^q\}$ in $\{0,1\}^n$, each with $m$ nonzero coordinates,
so that $q\leq {n\choose m}$.
Put $v^{q+1}=v^1$.
Further let $g\colon\mathbb{R}_+^m\to\mathbb{R}_+^m$ be a continuous order
preserving homogeneous map and assume that there exists $C>0$ such that
\[
g(z)_i\leq C(z_1\wedge z_2\wedge\ldots\wedge z_m)\mbox{\quad for all }
1\leq i\leq m\mbox{ and }z\in\mathbb{R}_+^m.
\]
Assume also that $g$ has a periodic point
$y$ with period $p$, where $1\leq p\leq {m\choose\lfloor
m/2\rfloor}$ and $\mathcal{O}(y;g)\subset
\mathrm{int}(\mathbb{R}_+^m)$. The existence of such a map $g$ and a
periodic point $y$ is guaranteed by taking a map $h$ as in 
 Lemma \ref{lem:6.1} and defining  
$g(z)_i=h(z)_i\wedge C(z_1\wedge\ldots\wedge z_m)$ for $1\leq i\leq m$, 
with $C$ large enough. 

For $1\leq k\leq q$  and $1\leq i\leq m$, we let $\nu(k,i)$ be the
index of the $i$th nonzero coordinate of $v^k$. Further for each
$x\in\mathbb{R}^n_+$, we let $x_{|v^k}$ be the vector in
$\mathbb{R}_+^m$ given by $(x_{|v^k})_i=x_{\nu(k,i)}$ for all 
$1\leq i\leq m$.  Subsequently,
we define $f\colon\mathbb{R}^n_+\to\mathbb{R}^n_+$ in the
following manner:
\begin{equation}
\label{eq:6.1} f(x)_i=\bigvee_{(k,r)\colon \nu(k+1,r)=i}
g(x_{|v^k})_r\mbox{\quad for each }1\leq i\leq n\mbox{ and }
x\in\mathbb{R}^n_+.
\end{equation}
It easy to see that $f\colon\mathbb{R}^n_+\to\mathbb{R}^n_+$ is a
continuous order preserving homogeneous map. Furthermore it has a
periodic point with period $\mathrm{lcm}\,(p,q)$. Indeed, for
$0\leq a\leq p-1$ and $1\leq b\leq q$ let
$y^{a,b}\in\mathbb{R}^n_+$ be given by 
\[
y^{a,b}_i=\left\{ \begin{array}{ll}
                   0      & \mbox{if }v^b_i=0\\
                 g^a(y)_r & \mbox{if }i=\nu(b,r).
                  \end{array}\right.
\]
As $\{g^k(y)\colon 0\leq k<p\}\subset\mathrm{int}(\mathbb{R}_+^m)$, it is
evident that
$y^{a,b}=y^{c,d}$ if and only if $a=c$ and $b=d$, so that they are all
distinct.
To complete the proof we now show that $f(y^{a,b})=y^{a+1,b+1}$, where the
indices $a$ and $b$ are counted modulo $p$ and modulo $q$, respectively.
As $g(z)_i\leq C(z_1\wedge z_2\wedge \ldots\wedge z_m)$ for each
$1\leq i\leq m$ and $z\in\mathbb{R}_+^m$, we have that
\[
g(y^{a,b}_{|v^k}) = \left\{ \begin{array}{ll}
                                          0      &\mbox{if }k\neq b\\
                                       g(g^a(y)) &\mbox{if }k=b.
                  \end{array}\right.
\]
Therefore
\[
f(y^{a,b})_i= \left\{ \begin{array}{ll}
                   0      & \mbox{if }v^{b+1}_i=0 \\
                 g^{a+1}(y)_r & \mbox{if }i=\nu(b+1,r),
                  \end{array} \right.
\]
for $1\leq i\leq n$.
Thus, $f(y^{a,b})=y^{a+1,b+1}$ and hence $y^{0,1}$ is a periodic point
of $f$ with period $\mathrm{lcm}\,(p,q)$.
\end{proof}
To illustrate the construction in the proof of Theorem \ref{thm:2.2}, we
consider the following example.
Let $m=2$, $n=3$, $p=2$, and $q=3$.
Put
\[
v^1=\left(\begin{array}{c} 1\\1\\0\end{array}\right),\,
v^2=\left(\begin{array}{c} 1\\0\\1\end{array}\right),\mbox{ and }
v^3=\left(\begin{array}{c} 0\\1\\1\end{array}\right).
\]
Further let $g\colon\mathbb{R}_+^2\to\mathbb{R}_+^2$ be given by
\[
g\left(\begin{array}{c}z_1\\z_2\end{array}\right)
=\left(\begin{array}{c}3z_1\wedge z_2 \\ z_1 \wedge
3z_2\end{array}\right),
\]
and take $y=(1,2)$.
It is easy to see that $y$ is a periodic point of $g$ with period 2.
The map $f\colon\mathbb{R}_+^3\to\mathbb{R}_+^3$ defined in (\ref{eq:6.1}) 
is then given by
\[
f\left(\begin{array}{c}x_1\\x_2\\x_3\end{array}\right)
= \left(\begin{array}{c} g(x_1,x_2)_1 \vee g(x_2,x_3)_1 \\ g(x_1,x_3)_1 \vee
g(x_2,x_3)_2 \\ g(x_1,x_2)_2 \vee g(x_1,x_3)_2 \end{array}\right)
=
\left(\begin{array}{c}
(3x_1\wedge x_2) \vee (3x_2\wedge x_3) \\
(3x_1\wedge x_3) \vee (x_2\wedge 3x_3) \\
(x_1\wedge 3x_2) \vee (x_1\wedge 3x_3) \end{array}\right).
\]
Now it is easy to verify that
\[
y^{0,1}= \left(\begin{array}{c}1\\2\\0\end{array}\right),\,
y^{1,2}=\left(\begin{array}{c}2\\0\\1 \end{array}\right),\,
y^{0,3}=\left(\begin{array}{c}0\\1\\2\end{array}\right),\,
\]
\[
y^{1,1}= \left(\begin{array}{c}2\\1\\0\end{array}\right),\,
y^{0,2}=\left(\begin{array}{c}1\\0\\2 \end{array}\right),\,
y^{1,3}=\left(\begin{array}{c}0\\2\\1\end{array}\right),\,
\]
is a periodic orbit of $f$ with period $\mathrm{lcm}\,(2,3)=6$.

We would also like to point out that if we take
$g(z)_i=z_1\wedge z_2\wedge\ldots\wedge z_m$ for all $1\leq i\leq m$ 
in the proof, then we recover the construction of Gunawardena and Sparrow. 
In particular, the maps $f$ and $g$ are so-called min-max maps. 

It follows directly from Theorem \ref{thm:2.2} that $\alpha_N$ given in
(\ref{eq:2.2})
is a lower bound for the maximal period of periodic points of continuous
order preserving subhomogeneous maps $f\colon K\to K$, where $K$ is a
polyhedral cone with $N$ facets in a finite dimensional vector space $X$.
By using the prime number theorem we now show that $\alpha_N$ has the same
asymptotics as the upper bound $\beta_N$ given in (\ref{eq:2.1}).
From Lemma \ref{lem:4.2} and equation (\ref{eq:4.2}) it follows that
\begin{equation}
\label{eq:6.2}
\beta_N=
\max_{\lfloor N/2\rfloor\leq m\leq N}{N\choose m}{m\choose \lfloor m/2\rfloor}=
\max_{1\leq m\leq N} {N\choose m}{m\choose \lfloor m/2\rfloor}.
\end{equation}
For a given $N$ let $m^*$ be the $m$ that attains the maximum in the
right-hand side of (\ref{eq:6.2}).
From the proof of Lemma \ref{lem:4.2} we know that
$m^*=\lfloor \frac{N+1}{3}\rfloor+\lfloor \frac{N+2}{3}\rfloor$.
Now for each $k\geq 1$ let $\rho(k)$ be the largest prime not exceeding $k$.
It then follows from the prime number theorem that
\begin{equation}
\label{eq:6.3}
\lim_{k\to\infty} \frac{\rho(k)}{k}=1.
\end{equation}
Indeed, let $p_N$ denote the $N$th prime and let $\pi(k)$ be the number of
primes not exceeding $k$.
Then $\rho(k)=p_{\pi(k)}$ for each $k\geq 1$.
It is known that
\[
\lim_{N\to\infty} \frac{p_N}{N\log N}=1\mbox{\quad and \quad }
\lim_{k\to\infty} \frac{\pi(k)\log \pi(k)}{k}=1
\]
are equivalent to the prime number theorem (see \cite[p.80]{Ap}).
Thus, the prime number theorem implies that
\[
\lim_{k\to\infty}\frac{\rho(k)}{k}=
\lim_{k\to\infty}\frac{p_{\pi(k)}}{\pi(k)\log \pi(k)}\cdot
\frac{\pi(k)\log \pi(k)}{k}=1.
\]
We now observe that
\[
\alpha_N\geq \rho({N\choose m^*}){m^*\choose \lfloor
m^*/2\rfloor},
\]
if $\rho({N\choose m^*})$ and ${m^*\choose \lfloor m^*/2\rfloor}$
are coprime. As $m^*=\lfloor \frac{N+1}{3}\rfloor+\lfloor
\frac{N+2}{3}\rfloor$, we can use Stirling's formula to show that
there exists $M\geq 1$ such that
\[
2{m^*\choose \lfloor m^*/2\rfloor}\leq {N\choose m^*} \mbox{\quad
for all }N\geq M \mbox{\quad and\quad } {N\choose m^*}\to
\infty,\mbox{ \quad as } N\to \infty.
\]
Therefore (\ref{eq:6.3}) implies that $\rho({N\choose m^*}) > 
{m^*\choose \lfloor m^*/2\rfloor}$ for all $N$
sufficiently large and hence they are coprime. Thus, we derive that
\[
\lim_{N\to\infty}\frac{\alpha_N}{\beta_N}\geq
\lim_{N\to\infty}\frac {{m^*\choose \lfloor
m^*/2\rfloor}\rho({N\choose m^*})} {{m^*\choose \lfloor
m^*/2\rfloor}{N\choose m^*}}=1.
\]
As $\alpha_N\leq\beta_N$ for each $N\geq 1$, we find that
$\lim_{n\to\infty} \alpha_N/\beta_N=1$.

We conclude the paper with some remarks. Given a polyhedral cone $K$, let
$\Gamma(K)$ be the set of integers $p\geq 1$ for which there exists
a continuous order preserving subhomogeneous map $f\colon K\to K$ that has a 
periodic point with period $p$. From Theorem \ref{thm:4.1} it follows that
$\Gamma(K)$ is a finite set. In fact, Theorem \ref{thm:4.1} implies that if 
$K$ has $N$ facets, then 
$\Gamma(K)\subset B(N)$, where $B(N)$ is the set of  
$p\geq 1$ for which there exist integers $q_1$ and $q_2$ such that
$p=q_1q_2$, $1\leq q_1\leq {N\choose m}$, and
$1\leq q_2\leq {m\choose \lfloor m/2\rfloor}$ for some $1\leq m\leq N$.
In particular, it follows that
$\Gamma(\mathbb{R}^3_+)\subset\{1,2,3,4,6\}$, so that 5 is not in 
$\Gamma(\mathbb{R}^3_+)$.
By Theorem \ref{thm:2.2} we know that $\Gamma(\mathbb{R}^N_+)\supset A(N)$,
where $A(N)$ is the set of $p\geq 1$ for which there
exist $1\leq m\leq N$, $1\leq q_1\leq {N\choose m}$, and $1\leq
q_2\leq {m\choose \lfloor m/2\rfloor}$ such that
$p=\mathrm{lcm}(q_1,q_2)$. For instance, 
$\Gamma(\mathbb{R}^3_+)\supset A(3)=\{1,2,3,6\}\neq B(3)$. 
Thus, for each $N\geq 1$ we have the following inclusions:
\[
A(N)\subset \Gamma(\mathbb{R}^N_+)\subset B(N).
\]
Knowing these inclusions it is natural to ask if there exists a
characterization of $\Gamma(\mathbb{R}^N_+)$ in terms of arithmetical and (or)
combinatorial constraints. In particular, one might wonder if
$\Gamma(\mathbb{R}^N_+)=A(N)$ for all $N\geq 1$, or, 
if $\Gamma(\mathbb{R}^N_+)=B(N)$ for all $N\geq 1$. 
This question is investigated by Bas Lemmens and Colin Sparrow in a 
forthcoming paper.

\end{document}